\documentclass[11pt]{article}
\usepackage{amsmath}
\usepackage{amsfonts}
\usepackage{amssymb}
\usepackage{amsthm}
\newtheorem*{theorem}{Theorem}
\newtheorem{lemma}{Lemma}
\title{Representation of One as the Sum of Unit Fractions}
\author{Yuya Dan\\Matsuyama University\\dan@cc.matsuyama-u.ac.jp}
\begin{document}
\maketitle

\begin{abstract}
One is expressed as the sum of the reciprocals of a certain set of integers.
We give an elegant proof to the fact applying the polynomial theorem and basic calculus.
\end{abstract}

\section{Introduction}
Let us consider the representation of one as the sum of unit fractions.
For examples,
we can take $2$, $3$ and $6$ for
\begin{equation}
	\frac{1}{2} + \frac{1}{3} + \frac{1}{6} = 1,
\end{equation}
and $3$, $4$, $4$, $8$ and $24$ for
\begin{equation}
	\frac{1}{3} + \frac{1}{4} + \frac{1}{4} + \frac{1}{8} + \frac{1}{24} = 1.
\end{equation}

It is well known that any positive rational number can be written as the sum of unit fractions.
In the paper,
we give a part of solutions to the Diophantine equation
\begin{equation}
	\frac{1}{x_1} + \frac{1}{x_2} + \cdots + \frac{1}{x_n} = 1,
\end{equation}
where the $x_j$ are not necessarily distinct integer for $j = 1, 2, \dots, n$.

To explain our result for $n = 6$,
we find all possible combinations of
\begin{equation}
	\alpha = ( \alpha_1, \alpha_2, \alpha_3, \alpha_4, \alpha_5, \alpha_6 ) \in \mathbb{N}^6
\end{equation}
such that
\begin{equation}
	\alpha_1 + 2 \alpha_2 + 3 \alpha_3 + 4 \alpha_4 + 5 \alpha_5 + 6 \alpha_6 = 6,
\end{equation}
where $\alpha_j \in \mathbb{N} = \{ 0, 1, 2, \dots \}$ for $j = 1$ to $6$. 
Next,
we take the quantities
\begin{equation}
	\prod_{j=1}^{6} \alpha_j! j^{\alpha_j} = \alpha_1! 1^{\alpha_1} \cdot \alpha_2! 2^{\alpha_2} \cdot \cdots \cdot \alpha_6! 6^{\alpha_6}
\end{equation}
for each possible $\alpha$.
Then,
we can calculate the sum of reciprocals of the above quantities
\begin{equation}
	\frac{1}{720} + \frac{1}{48} + \frac{1}{18} + \frac{1}{16} + \frac{1}{8} + \frac{1}{6} + \frac{1}{5} + \frac{1}{48} + \frac{1}{8} + \frac{1}{18} + \frac{1}{6}
\end{equation}
which is equal to 1.

\begin{table}
\caption{Possible combinations of $\alpha$ for $n = 6$}
\begin{center}
\begin{tabular}{ccccccl}
	\hline
	$\alpha_1$ & $\alpha_2$ & $\alpha_3$ & $\alpha_4$ & $\alpha_5$ & $\alpha_6$ & denominator \\
	\hline
	6 & 0 & 0 & 0 & 0 & 0 & $6! = 720$ \\
	4 & 1 & 0 & 0 & 0 & 0 & $4! \cdot 2 = 48$ \\
	3 & 0 & 1 & 0 & 0 & 0 & $3! \cdot 3 = 18$ \\
	2 & 2 & 0 & 0 & 0 & 0 & $2! \cdot 2! \cdot 2^2 = 16$ \\
	2 & 0 & 0 & 1 & 0 & 0 & $2! \cdot 4 = 8$ \\
	1 & 1 & 1 & 0 & 0 & 0 & $2 \cdot 3 = 6$ \\
	1 & 0 & 0 & 0 & 1 & 0 & $5$ \\
	0 & 3 & 0 & 0 & 0 & 0 & $3! \cdot 2^3 = 48$ \\
	0 & 1 & 0 & 1 & 0 & 0 & $2 \cdot 4 = 8$ \\
	0 & 0 & 2 & 0 & 0 & 0 & $2! \cdot 3^2 = 18$ \\
	0 & 0 & 0 & 0 & 0 & 1 & $6$ \\
	\hline
\end{tabular}
\end{center}
\end{table}

The example is generalized to our main result:
\begin{theorem}
For any positive integer $n$,
\begin{equation}
	\sum_{\alpha \in S_n} \prod_{j=1}^n \frac{1}{\alpha_j! j^{\alpha_j}} = 1,
\end{equation}
where the summation over $S_n$ runs through all possible $\alpha = ( \alpha_1, \alpha_2, \dots, \alpha_n )$ in $\mathbb{N}^n$ such that
\begin{equation}
	\sum_{j=1}^n j \alpha_j = n.
\end{equation}
\end{theorem}

\section{Preliminaries}
\begin{lemma}[Polynomial theorem]
\label{PolynomialTheorem}
Let n and m be positive integers.
For any $x = ( x_1, x_2, \dots, x_n )$ in $\mathbb{R}^n$,
\begin{equation}
	( x_1 + x_2 + \cdots + x_n )^m = \sum_{\lvert \alpha \rvert = m} \frac{m!}{\alpha!} x^{\alpha}
\end{equation} 
where $\alpha! = \alpha_1! \alpha_2! \cdots \alpha_n!$ for $\alpha = ( \alpha_1, \alpha_2, \dots, \alpha_n )$ in $\mathbb{N}^n$ and the summation runs through all possible $\alpha$ in $\mathbb{N}^n$ such that $\lvert \alpha \rvert = \alpha_1 + \alpha_2 + \cdots + \alpha_n = m$.
\end{lemma}

\begin{proof}
Each coefficient of
\begin{equation}
	x^{\alpha} = \prod_{j=1}^n x_j^{\alpha_j}
\end{equation}
in the right-hand side for some $\alpha$ in $\mathbb{N}^n$ with $\lvert \alpha \rvert = m$ is equal to the number of combinations of the products among $x_1, x_2, \dots, x_n$.
\end{proof}

\begin{lemma}
\label{Coefficient}
Given a polynomial
\begin{equation}
	f( x ) = \sum_{j = 0}^n a_j x^j.
\end{equation}
Then,
the $j^\text{th}$ coefficient of $f( x )$ can be expressed by
\begin{equation}
	a_j = \frac{1}{j!} f^{(j)}( 0 ),
\end{equation}
where $f^{(j)}$ stands for the $j^\text{th}$ derivative.
\end{lemma}

\begin{proof}
$f$ is infinitely differentiable,
since the $j$-times differential function of $x^k$ is
\begin{equation}
	\left( \frac{d}{dx} \right)^j x^k = \frac{k!}{( k - j )!} x^{k - j}
\end{equation}
if $j \leq k $,
and
\begin{equation}
	\left( \frac{d}{dx} \right)^j x^k = 0
\end{equation}
if $j > k$.
We have
\begin{equation}
	f^{(j)}( x )
	= \sum_{k=0}^n a_k \left( \frac{d}{dx} \right)^j x^k
	= \sum_{k=j}^n \frac{k!}{( k - j )!} a_k x^{k - j}
\end{equation}
for any $j$ between $0$ and $n$,
then
\begin{equation}
	f^{(j)}( 0 ) = j! a_j,
\end{equation}
which implies the conclusion of the lemma
\begin{equation}
	a_j = \frac{1}{j!} f^{(j)}( 0 ).
\end{equation}
\end{proof}

\begin{lemma}
\label{Coefficient2}
Let $n$ be a positive integer.
We put
\begin{equation}
	g( x ) = \left( x + \frac{1}{2} x^2 + \cdots + \frac{1}{n} x^n \right)^n,
\end{equation}
then $g^{(n)}( 0 ) = n!$.
\end{lemma}

\begin{proof}
Put
\begin{equation}
	g( x ) = x^n \left( 1 + \frac{1}{2} x + \cdots + \frac{1}{n} x^{n - 1} \right)^n
	= x^n h( x ),
\end{equation}
where
\begin{equation}
	h( x ) = \left( 1 + \frac{1}{2} x + \cdots + \frac{1}{n} x^{n - 1} \right)^n.
\end{equation}
Leibniz rule implies
\begin{equation}
\begin{array}{ccl}
	g^{(n)}( x ) &=& \displaystyle\sum_{j = 0}^n \frac{n!}{( n - j )! j!} \left( \frac{d}{dx} \right)^{n - j} x^n \cdot \left( \frac{d}{dx} \right)^{j} h( x )\\
	&=& \displaystyle\sum_{j = 0}^n \frac{n! n!}{( n - j )! j! j!} x^j \cdot \left( \frac{d}{dx} \right)^{j} h( x ).
\end{array}
\end{equation}
Therefore,
we obtain
\begin{equation}
	g^{(n)}( 0 ) = n! h( 0 ) = n!.
\end{equation}
\end{proof}

\section{Proof of the main result}

We begin with the relation
\begin{equation}
\begin{array}{ccl}
	\displaystyle\left( x_1 + \frac{1}{2} x_2^2 + \cdots + \frac{1}{n} x_n^n\right)^m
	&=& \displaystyle\sum_{\lvert \alpha \rvert = m} \frac{m!}{\alpha!} x_1^{\alpha_1} \cdot \left( \frac{1}{2} x_2^2 \right)^{\alpha_2} \cdot \cdots \cdot \left( \frac{1}{n} x_n^n \right)^{\alpha_n}\\
	&=& m! \displaystyle\sum_{\lvert \alpha \rvert = m} \prod_{j=1}^n \frac{1}{\alpha_j ! j^{\alpha_j}} x_j^{j \alpha_j}.
\end{array}
\end{equation}
by Lemma \ref{PolynomialTheorem}.
Putting $m = n$ and $x_j = t$ for $j = 1, 2, \dots, n$ implies
\begin{equation}
	\left( t + \frac{1}{2} t^2 + \cdots + \frac{1}{n} t^n \right)^n
	= n! \sum_{\lvert \alpha \rvert = n} \prod_{j=1}^n \frac{1}{\alpha_j ! j^{\alpha_j}} t^{j \alpha_j}.
\end{equation}
Compare to the coefficients of $t^{n}$ in both side of the identity.
We obtain $n!$ from the left-hand side of the identity by Lemma \ref{Coefficient} and Lemma \ref{Coefficient2}.
On the other hand,
the coefficient of $t^{n}$ in the right-hand side is the sum of all terms with $t^{n}$,
which is written by
\begin{equation}
	\sum_{\alpha \in S_n} \prod_{j=1}^n \frac{1}{\alpha_j ! j^{\alpha_j}},
\end{equation}
where the summation runs through all possible $\alpha$ in $S_n$ defined by
\begin{equation}
	S_n = \{ \alpha \in \mathbb{N}^n; \alpha_1 + 2 \alpha_2 + \cdots + n \alpha_n = n \}.
\end{equation}
Hence,
we obtain
\begin{equation}
	\sum_{\alpha \in S_n} \prod_{j=1}^n \frac{1}{\alpha_j ! j^{\alpha_j}} = 1
\end{equation}
which completes the proof of our main result.

\section{Concluding Remarks}
In this paper,
a part of reciprocal bases of one is investigated from analytic point of view.
In particular,
the polynomial theorem and the multi-index analysis play important role in the proof.
Although these are not all of solutions to the Diophantine equation,
one is presented as the sum of the reciprocal numbers of a certain set of integers.

\end{document}